\documentclass{birkjour}
\usepackage{amsmath}
\usepackage{amssymb}
\usepackage{amsthm}
\usepackage[all]{xy}
\usepackage{enumitem}
\usepackage[T1]{fontenc}

 \newtheorem{thm}{Theorem}[section]
 
 \newtheorem{lem}[thm]{Lemma}
 \newtheorem{prop}[thm]{Proposition}
 \theoremstyle{definition}
 \newtheorem{defn}[thm]{Definition}
 \theoremstyle{remark}
 \newtheorem{rem}[thm]{Remark}
 \newtheorem*{ex}{Examples}
 \numberwithin{equation}{section}

\begin{document}

\title[Hamiltonian dynamics on multicotangent bundles]{Notes on equivalent formulations of Hamiltonian dynamics on 
multicotangent bundles}


\author[Wagner]{Maxime Wagner}
\address{Institut \'Elie Cartan de Lorraine\br
Universit\'e de Lorraine et CNRS\br
57070 Metz, France}
\email{maxime.wagner@univ-lorraine.fr}

\author[Wurzbacher]{Tilmann Wurzbacher}
\address{Institut  \'Elie  Cartan de Lorraine\br
Universit\'e de Lorraine et CNRS\br
57070 Metz, France}
\email{tilmann.wurzbacher@univ-lorraine.fr}



\begin{abstract}
We show the equivalence of five different conditions on a classical field $\psi$ with values in a restricted multicotangent bundle to be a solution of the field equations, notably in terms of the Hamilton-Volterra equations, the principle of least action and several conditions based on the contraction of the multi-vector tangent to $\psi$ with canonical differential forms. Most prominently, we have equivalence to the ``dynamical Hamilton-de Donder-Weyl equation'', that can be vastly generalized to define Hamiltonian dynamics on multisymplectic manifolds, defined for sources of different dimensions.
\end{abstract}

\maketitle
\tableofcontents
\clearpage

\section{Introduction}\label{Introduction}

The main purpose of this text is to explain five rather different but equivalent conditions for a \textit{classical field} to be a solution of the field equations in the Hamiltonian approach, within the widely used standard framework of a fibre bundle 
$E\stackrel{\pi}{\rightarrow}\Sigma$ over $\Sigma$, the source of the fields. A classical field $\psi$ is here a section of the restricted multimoment bundle. 
This corresponds often to the Lagrangian approach, where a field is a section of $\pi$, via the reduction of order and a Legendre transformation applied to the Lagrangian side. Since Lagrangian functions or densities are not always sufficiently non-degenerate to allow for a bijective Legendre transformation, there is not a one-to-one correspondence between Lagrangian and Hamiltonian classical field theories. Having in mind \textit{Hamiltonian dynamics} on multisymplectic manifolds, generalizing dynamical systems on symplectic manifolds, we do not consider the Lagrangian approach in this article but concentrate entirely on the Hamiltonian approach. Furthermore we concentrate on $M(\pi)=\Lambda^n_2T^*E$,
the restricted multicotangent bundle associated to $\pi$ and its quotient $P(\pi)=\Lambda^n_2T^*E/\Lambda^n_1T^*E$, the
natural field theoretic analogues of cotangent bundles in time-dependent classical mechanics.
\\\\
Our second important goal is to generalize the considered field equation on $M(\pi)$ to a vast class of equations, we call ``dynamical Hamilton-de Donder-Weyl equation'', for maps from a $k$-dimensional source $\Sigma$ equipped with a co-volume to an $n$-plectic manifold $M$ with $1\leq k\leq n$. (Note that the existence of a co-volume on $\Sigma$ guarantees in the muticotangent case that we can work with a Hamiltonian function instead of a Hamiltonian section, cf. below. Since a general multisymplectic manifold is not necessarily equipped with any natural fibre bundle structure, a reformulation in terms of a Hamiltonian function on $M(\pi)$ is in fact needed for the above generalization).
\\\\
In Section 2 we explain the set-up of these multicotangent bundles and the basics of their geometry. In Section 3 we clarify the relation between Hamiltonian sections $h$ on $P(\pi)$ and Hamiltonian functions $H$ on $M(\pi)$, since it is indispensable for the generalization of the field equations from
multicotangent bundles to arbitrary multisymplectic manifolds.
\\\\
We then give, in Section 4, precise mathematical formulations of five different conditions for a classical field $\psi$ \big(or its quotient version $\widetilde{\psi}=\rho \circ \psi$, $\rho$ being the projection $M(\pi)\to P(\pi)$\big). Historically, the first formulation is certainly the principle of least action, duly generalized to classical field theories and going back, in mechanics, at least to Maupertuis and Hamilton. In the field theoretic context, Volterra gave already in 1890 a generalization of Hamilton's equation of mechanics, nowadays called Hamilton-Volterra equations.
\\\\
We pay special attention to the fact that in the variational formulation we have to work with a local version of criticality (coined ``extremality'') replacing the here ill-defined notion of (global) criticality.
\\\\
Furthermore, we characterize solutions via the vanishing of certain canonical differential forms on the image of $\psi$ or equivalently of their pull-back to the source manifold $\Sigma$. This yields the following condition, put forward by the ``Spanish school'' since the 1990's but already present in special cases in \cite{Gawedzki}\::

\begin{eqnarray*}
    \widetilde{\psi}^*\big(\iota_X(h^*\omega)\big)=0 \hskip0.3cm \text{ for all vector fields } X \text{ on }P(\pi).
\end{eqnarray*}
This equation translates to the condition that $\psi = h \circ \widetilde{\psi}$ is a vortex $n$-plane, generalizing the vortex lines of Arnold (cf. \cite{Arnold}, p. 235), and to the closely related following condition we call \textit{dynamical Hamilton-de Donder-Weyl equation}

\begin{eqnarray*}
    \iota_{(\psi_*)_x(\gamma_x)}\omega_{\psi(x)}=(dH)_{\psi(x)} \hskip0.4cm \forall x\in \Sigma,
\end{eqnarray*}
where $\gamma$ is a multi-vector field of degree $n$, dual to a fixed volume form on the source $\Sigma$.
\\\\
This last equation allows for a vast generalization to arbitrary $n$-plectic manifolds (i.e. manifolds equipped with a multisymplectic form of degree $n{+}1$) and to source manifolds $\Sigma$ of dimension $k$ with $1\leq k \leq n$. These equations, notably for $k=1$ and $k=n$ are, of course,  central to the promising emerging subject of \textit{Hamiltonian dynamics on multisymplectic manifolds}, whose genesis starting from classical field theory on multicotangent bundles we describe in this article.
\\\\
Certain results in this article are of a folkloric nature but their proofs are widely scattered in the literature, with varying degrees of rigour in the formulation of the conditions and in the proof of their equivalences. We give here complete definitions and proofs, and explain in detail the transition to fields with values in arbitrary multisymplectic manifolds \big(cf. notably Theorem \ref{section dynamical HDW} and Definition \ref{Final Definition} (iii)\big). On the other hand we apologize for not being exhaustive in our list of references, but we did not try to describe the highly interesting history of multisymplectic geometry.

\section{Multicotangent bundles and classical field theory}\label{Section : Multicotangent bundles and Classical Field Theory}

In this section we introduce the restricted multicotangent bundle $M(\pi)$ and its quotient $P(\pi)$, direct generalizations of $T^*(I\times Q)$ and $I\times T^*Q$, the phase spaces of time-dependent mechanics (compare \cite{Carinena&Crampin&Ibort}, \cite{Echeverrai&MunozLecanda&RomanRoy2000}, \cite{Forger&Paufler&Romer}, \cite{Goldschmidt&Sternberg}, \cite{Kijowski&Tulczyjew}, \cite{Ryvkin&Wurzbacher} for the spaces $M(\pi)$ and $P(\pi)$ and their relations to appropriate duals of jets bundles). We also introduce the highly useful \textit{standard coordinate systems} on them and show that they can locally be adapted to a given volume form on the base $\Sigma$, generalizing the interval $I$ of the mechanics.

\begin{defn}
    Let $\pi:E\rightarrow \Sigma$ be a smooth fibre bundle with typical fibre $Q$ of dimension $N$ and let $n=\dim \Sigma$. For $1\leq k\leq n+1$ we set 
    \begin{eqnarray*}
    \Lambda^n_kT^*E=\{\eta\in \Lambda^nT^*E\,|\, \eta \text{ lies over } y\in E \\
    \text{ and } \forall u_1,\ldots ,u_k \in \negthinspace  \negthinspace \negthinspace & \negthinspace\ker(\pi_*)_y,\,
    \iota_{u_k}\ldots\iota_{u_1}\eta=0\}.
    \end{eqnarray*}
\end{defn}

In the sequel we use the following notations :
\begin{align}\label{diagram}
    \xymatrix{
    \Lambda^n_2T^*E= M(\pi) \ar[d]_{\rho} \\
    \frac{\Lambda^n_2T^*E}{\Lambda^n_1T^*E}= P(\pi) \ar@/_2pc/[dd]_{\tau} \ar[d]_{\kappa} \\
     E \ar[d]_{\pi}^Q\\
     \Sigma.
     }
\end{align}

\begin{rem}
    The fibre bundle $\pi:E\rightarrow \Sigma$ is of course assumed to be locally trivial but it is not necessarily a vector bundle or a principal bundle. Thus we include, e.g. classical nonlinear sigma models here.
\end{rem}

\begin{defn}
    In the above situation, let $U\subset \Sigma$ a chart domain with coordinate map $x=(x^1,\ldots,x^n)$ together with a trivialization $E\big|_U\xrightarrow{\zeta_U} U\times Q$ and let $W\subset Q$ be a chart domain with coordinate map $q=(q^1,\ldots q^N)$. Then we call :
    \begin{enumerate}[label=(\roman*)]
        \item the coordinate chart $(x,\,q)\circ\zeta_U:V=\zeta_U^{-1}(U\times W)\rightarrow x(U)\times q(W)\subset \mathbb{R}^n\times \mathbb{R}^N$ \textit{standard} \big(in the sequel, we will write $(x,\,q)$ for $(x,\,q)\circ\zeta_U$\big)\;;\\
        \item a coordinate chart $(x^{\mu},\,q^a,\,p^{\mu}_a,\,p)$ on $(\kappa\circ\rho)^{-1}(V)\subset M(\pi)$ \textit{standard} if $(x,\,q)$ is standard on $V$ and for all $\eta\in (\kappa\circ\rho)^{-1}(V)$, we have $\eta=pd^nx+\sum\limits_{a=1}^N\sum\limits_{\mu=1}^np^{\mu}_adq^a\wedge d^{n-1}\widehat{x^{\mu}}$ where $d^{n-1}\widehat{x^{\mu}}=\iota_{\frac{\partial}{\partial x^{\mu}}}d^nx=(d^nx)(\frac{\partial}{\partial x^{\mu}},-)$\;;\\
        \item a coordinate system $(x^{\mu},\,q^a,\,p^{\mu}_a)$ on $\kappa^{-1}(V)\subset P(\pi)$, \textit{standard} if $(x,\,q)$ is standard on $V$ and for all $[\eta]\in \kappa^{-1}(V)$ with $\eta\in (\kappa\circ\rho)^{-1}(V)$, we have $[\eta]=\sum\limits_{a=1}^N\sum\limits_{\mu=1}^np^{\mu}_adq^a\wedge d^{n-1}\widehat{x^{\mu}}$ modulo $\mathbb{R}(d^nx)$.
    \end{enumerate}
\end{defn}

\begin{rem}
    It is obvious that such standard coordinate systems exist.
\end{rem}

\begin{defn}
    Suppose a volume form ${vol}^{\Sigma}$ is given on $\Sigma$. A coordinate chart is a \textit{standard coordinate chart adapted to} ${vol}^{\Sigma}$ if the chart is standard in the previous sense and $d^nx={vol}^{\Sigma}$.
\end{defn}

\begin{lem}
    Let $\Sigma$ a $n$-dimensional manifold,  $\hbox{vol}^{\Sigma}$ a volume form on $\Sigma$ and $x_0\in \Sigma$. Then there exists a neighborhood $U$ of $x_0$ and a coordinate chart $(x^1,\ldots,x^n)$ on $U$ such that ${vol}^{\Sigma}_{|U}=d^nx$.
\end{lem}

\begin{proof}
    \raggedright Let $(\Tilde{x}^1,...,\Tilde{x}^n)$ be a coordinate system near $x_0$. There exists $f$ with no zeros near $x_0$ such that ${vol}^{\Sigma}=f(\widetilde{x})d\Tilde{x}^1\wedge \ldots \wedge d\Tilde{x}^n$. Let $F(\Tilde{x})$ be a function, such that near $x_0$, $\frac{\partial F}{\partial \Tilde{x}^1}(\Tilde{x})=f(\Tilde{x})$, i.e. 
    \begin{eqnarray*}
        F(\Tilde{x})=\int f(\Tilde{x}_1,\,\Tilde{x}_2,\ldots,\,\Tilde{x}_n)\;d\Tilde{x}_1
    \end{eqnarray*} 
    is a primitive of $f$ with respect to $\Tilde{x}_1$. Defining $x^1=F(\Tilde{x})$, $x^2=\Tilde{x}^2$,\,\ldots, $x^n=\Tilde{x}^n$, it follows that ${vol}^{\Sigma}= 
    dF(\Tilde{x})\wedge d\Tilde{x}^2\wedge \ldots \wedge d\Tilde{x}^n=dx^1\wedge \ldots\wedge dx^n$ near $x_0$.
\end{proof}

\section{Hamiltonian sections versus Hamiltonian functions}\label{Section : Hamiltonian sections versus Hamiltonian functions}

In the absence of a fixed volume form on the base manifold $\Sigma$ (i.e. the source of the classical fields) a Lagrangian is a map from jets to densities. This is reflected in the multicotangent bundle picture by being forced to consider a Hamiltonian section $h$ of the bundle $M(\pi)\xrightarrow{\rho} P(\pi)$. As in the case of Lagrangian theory we can go to a scalar \textit{Hamiltonian function} on $M(\pi)$ in case we have fixed a volume ${vol}^{\Sigma}$ on $\Sigma$. This transition, only hinted at in the existing literature, is here worked out in all details. Notably, we characterize which Hamiltonian functions on $M(\pi)$ arise from Hamiltonians sections. (The notion of a Hamiltonian section is, i.e. developed in \cite{Carinena&Crampin&Ibort} and \cite{Echeverria&MunozLecanda&RomanRoy1999}, compare also \cite{Ryvkin&Wurzbacher}.)
\\\\
Let us first observe the following two obvious lemmata.

\begin{lem}
    Let $V$ a vector space and $W\subset V$ a linear subspace. Then
    \begin{eqnarray*}
        \rho : V & \rightarrow & V/ W\\
        v & \mapsto & [v]
    \end{eqnarray*}
    is an affine bundle modelled on the (trivial) vector bundle $V/ W \times W \rightarrow V / W$.
\end{lem}

\begin{lem}\label{fibre Lem}
    Let $\pi^E:E\rightarrow B$ be a vector bundle with typical fibre $V$ and $\pi^F : F \rightarrow B$ a vector sub-bundle of $E$ with typical fibre $W\subset V$. Then
    \begin{enumerate}[label=(\roman*)]
        \item 
        \begin{eqnarray*}
            \rho : E & \rightarrow & E/ F\\
            v_b & \mapsto & [v_b]
        \end{eqnarray*}
        defined by $\forall b\in B,\, \forall v_b\in E_b=\big(\pi^E\big)^{-1}(b)$
        is an affine bundle modelled on $(\pi^{E/ F})^*F$ with respect to the action\::
        \raggedright\begin{align*}
            \Theta : (\pi^{E/F})^*F\times_{E/F} E & \rightarrow E\\
            \big( ([v_b],\,u_b), v_b\big) & \mapsto v_b+u_b
         \end{align*}
        where $\pi^{E/ F} : E/ F \rightarrow B$ is the induced bundle projection\;;\\
        \item if $h$ is a smooth section of $\rho$, then we have the following diagram 
        $$\xymatrix{
        (\pi^{E/F})^*F \ar[rr]^{\Xi} \ar[rd] & & E \ar[ld]^\rho\\
        & E/F \ar@/_1pc/[ur]_h}$$ 
        with $\Xi$ defined by
        \begin{eqnarray*}
            \forall b\in B,\,  \forall [v_b]\in (E/ F)_b,\, \Xi_{[v_b]} (u_b) = \Theta \big(([v_b],\,u_b),\,h([v_b])\big) =h([v_b])+u_b
        \end{eqnarray*}
        being an isomorphism of affine bundle.
    \end{enumerate}
\end{lem}

\begin{prop}\label{Def chi}
    In the situation of Diagram \ref{diagram},
    \begin{enumerate}[label=(\roman*)]
        \raggedright\item there is a vector bundle isomorphism between $\Lambda^n_1T^*E$ and $\pi^*(\Lambda^nT^*\Sigma)$\;;\\
        \item the fibre bundle $\rho : M(\pi)\rightarrow P(\pi)$ is an affine bundle modelled on the vector bundle $\kappa^*(\Lambda^n_1T^*E)\rightarrow P(\pi)$\;;\\
        \item if $h$ is a section of $\rho : M(\pi)\rightarrow P(\pi)$, then $\rho$ is isomorphic to $\kappa^*(\Lambda^n_1T^*E)\rightarrow P(\pi)$ as an affine bundle\;;\\
        \item if there is fixed a volume form ${vol}^{\Sigma}$ on $\Sigma$, then $\Lambda^nT^*\Sigma$ is canonically trivializable\::
        $$\xymatrix{
        \Sigma\times\mathbb{R} \ar[rr]^{\phi} \ar[rd] & & \Lambda^nT^*\Sigma \ar[ld]\\
        & \Sigma
        }$$
        and, in this case, $\pi^*(\Lambda^nT^*\Sigma)\rightarrow E$ and $\tau^*(\Lambda^nT^*\Sigma)\rightarrow P(\pi)$ are also canonically trivializable\;;\\
        \item if $h$ is a section of $\rho$ and there is fixed a volume form on $\Sigma$ noted ${vol}^{\Sigma}$, then $\rho$ is canonically trivializable, i.e. we have the following diagram\::
        $$\xymatrix{
        P(\pi)\times \mathbb{R} \ar@/^2pc/[rr]^{\chi} \ar[rd] \ar[r]^{\cong} & \tau^*(\Lambda^nT^*\Sigma) \ar[d] \ar[r]^{\cong} & M(\pi) \ar[ld]_{\rho}\\
        & P(\pi) \ar@/_1pc/[ur]_h
        }$$
        where $\chi(z,\,u)=h(z)+u\cdot \tau^*({vol}^{\Sigma})_z$. (The horizontal isomorphisms are isomorphisms of affine bundles over $P(\pi)$.)
    \end{enumerate}
\end{prop}

\begin{proof}
    \begin{enumerate}[label=(\roman*)]
        \raggedright\item Let $y\in E$, and define 
        \begin{eqnarray*}
            \beta : \pi^*(\Lambda^nT^*\Sigma) \rightarrow \Lambda^n_1T^*E
        \end{eqnarray*}
        in the following way\:: for $v_1,\ldots,v_n\in T_yE$ we put 
        \begin{eqnarray*}
            \beta(\pi^*vol^{\Sigma})_y(v_1,\ldots,v_n)=vol^{\Sigma}_{\pi(y)}\big((\pi_*)_y(v_1),\ldots,(\pi_*)_y(v_n)\big).
        \end{eqnarray*}
        Then $\beta$ is injective and since rank $\pi^*(\Lambda^nT^*\Sigma)=1=$ rank $\Lambda^n_1T^*E$, this is a vector bundle isomorphism.\\
        \item Apply Lemma \ref{fibre Lem}.\\
        \item Apply Lemma \ref{fibre Lem}.\\
        \item The map $\phi$ defined by
        \begin{eqnarray*}
            \phi : \Sigma\times \mathbb{R} & \rightarrow & \Lambda^nT^*\Sigma\\
            (x,\,u) & \mapsto &u\cdot({vol}^{\Sigma})_x.
        \end{eqnarray*}
        is a trivialization of the fibre bundle $\Lambda^nT^*\Sigma\rightarrow \Sigma$.\\
        Obviously $\pi^*(\Lambda^nT^*\Sigma)\rightarrow E$, respectively $\tau^*(\Lambda^nT^*\Sigma)\rightarrow P(\pi)$, are then also trivializable via the maps
        \begin{eqnarray*}
            \phi^E : E\times \mathbb{R} & \rightarrow & \pi^*(\Lambda^nT^*\Sigma)\\
            (y,\,u) & \mapsto &u\cdot\big(\pi^*({vol}^{\Sigma})\big)_y,
        \end{eqnarray*}
        respectively,
        \begin{eqnarray*}
            \phi^{P(\pi)} : P(\pi)\times \mathbb{R} & \rightarrow & \tau^*(\Lambda^nT^*\Sigma)\\
            (z,\,u) & \mapsto &u\cdot\big(\tau^*({vol}^{\Sigma})\big)_z.
        \end{eqnarray*}
        \item As $h$ is a section of $\rho$, there is a canonical affine bundle isomorphism between $M(\pi)\rightarrow P(\pi)$ and $\tau^*(\Lambda^nT^*\Sigma)\rightarrow P(\pi)$. Then, the map
        \begin{eqnarray*}
            \chi : P(\pi)\times \mathbb{R} & \rightarrow &M(\pi)\\
            (z,\,u) & \mapsto &h(z)+u\cdot (\tau^*{vol}^{\Sigma})_z
        \end{eqnarray*}
        is an affine bundle isomorphism between $P(\pi)\times \mathbb{R}$ and $M(\pi)$.
    \end{enumerate}
\end{proof}

\begin{defn}
    We say that $h$ is a \textit{Hamiltonian section} if $h$ is a section of $\rho : M(\pi) \rightarrow P(\pi)$.
\end{defn}

\begin{prop}\label{Proposition H=h+p}
    \raggedright Let ${vol}^{\Sigma}$ be a volume form on $\Sigma$, $(x^{\mu},\,q^a,\,p^{\mu}_a,\,p)$ a standard coordinate chart adapted to ${vol}^{\Sigma}$ on $M(\pi)$ and $h$ a Hamiltonian section. We can define a Hamiltonian function $H\::M(\pi)\rightarrow \mathbb{R}$ via $H=\text{pr}_{\mathbb{R}}\circ \chi^{-1}$, where $\chi$ is the map in (v) of Proposition \ref{Def chi}.\\
    In coordinates, if $h(x^{\mu},\,q^a,\,p^{\mu}_a)=\big(x^{\mu},\,q^a,\,p_a^{\mu},-\mathcal{H}(x^{\mu},\,q^a,\,p_a^{\mu})\big)$, with a locally defined function $\mathcal{H}\:: P(\pi)\rightarrow \mathbb{R}$, we have $H(x^{\mu},\,q^a,\,p^{\mu}_a,\,p)=\mathcal{H}(x^{\mu},\,q^a,\,p^{\mu}_a)+\nolinebreak[4] p$.
\end{prop}

\begin{proof}
    We recall $\chi(z,\,u)=h(z)+u(\tau^*{vol}^{\Sigma})_z$ and in standard coordinates adapted to the volume form, ${vol}^{\Sigma}=d^nx$. Locally, we can write 
    \begin{eqnarray*}
        h : P(\pi) & \longrightarrow & M(\pi)\\
        \Big[\sum_{\mu,a}p^{\mu}_adq^a\wedge d^{n-1}\widehat{x^{\mu}}\Big] & \longmapsto & -\mathcal{H}(x^{\mu},\,q^a,\,p^{\mu}_a)d^nx+\sum_{\mu,a}p^{\mu}_adq^a\wedge d^{n-1}\widehat{x^{\mu}}.
    \end{eqnarray*}
    Then,  $\forall z\in P(\pi)$ and $u\in \mathbb{R}$,
    \begin{eqnarray*}
        \chi(z,\,u) & = & -\mathcal{H}(z)(d^nx)_z+\sum\limits_{\mu,a}p^{\mu}_a(dq^a\wedge d^{n-1}\widehat{x^{\mu}})_z+u(\tau^*{vol}^{\Sigma})_z\\
        & = & (-\mathcal{H}(z)+u)(d^nx)_z+\sum\limits_{\mu,a}p^{\mu}_a(dq^a\wedge d^{n-1}\widehat{x^{\mu}})_z.
    \end{eqnarray*}
    (Note that, by a slight abuse of language, $\tau^*{vol}^{\Sigma}$ is written as $d^nx$ here.)\\
    We observe that $p=-\mathcal{H}(z)+u$, so, $\forall z\in P(\pi)$, $u=\mathcal{H}(z)+p$. Thus, in local coordinates, $H$ is given as claimed.
\end{proof}

\begin{rem}
    The (locally defined) function $\mathcal{H}$ is also called the \textit{(local) Hamilton-Volterra function}.
\end{rem}

\begin{defn}
    Let $\eta \in M(\pi)=\Lambda^n_2T^*E$ and $f : M(\pi)\rightarrow \mathbb{R}$ be a smooth function. We define the vector field $Z$ on $M(\pi)$ by
    \begin{eqnarray*}
        Z_{\eta}(f)=\frac{d}{dt}\Big|_{t=0}f\Big(\eta+t\big(\pi^*{vol}^{\Sigma}\big)_{\kappa\circ\rho(\eta)}\Big).
    \end{eqnarray*}
\end{defn}

\begin{lem}\label{Flot Z}
    \raggedright Let $H\in\mathcal{C}^{\infty}(M(\pi),\,\mathbb{R})$ be a Hamiltonian function such that $Z(H)=1$ and note by $(\varphi_u^Z)_{u\in\mathbb{R}}$ the flow associated to $Z$. Then $\forall \eta \in M(\pi), \forall u\in \mathbb{R}$, we have $H\big(\varphi_u^Z(\eta)\big)=H(\eta)+u$.
\end{lem}

\begin{proof}
    \raggedright The flow associated to $Z$ is obviously given by $\varphi_u^Z(\eta)=\eta+u(\pi^*{vol}^{\Sigma})_{\kappa\circ\rho(\eta)}$. For a fixed $\eta$,  
    \begin{eqnarray*}
        H\big(\varphi_u^Z(\eta)\big)-H(\eta) & = & H\big(\varphi_u^Z(\eta)\big)-H\big(\varphi_0^Z(\eta)\big)\\
        & = & \int_0^u \frac{d}{dt}H\big(\varphi_t^Z(\eta)\big)dt\\
        & = & \int_0^u Z_{\varphi_t^Z(\eta)}(H)dt=\int_0^u 1dt=u.
    \end{eqnarray*}
\end{proof}

\begin{prop}
    Suppose that $\Sigma$ is orientable. Then we have
    \begin{enumerate}[label=(\roman*)]
        \raggedright \item given a Hamiltonian section $h$, the Hamiltonian function $H\:: M(\pi)\rightarrow \mathbb{R}$ defined by $H=\text{pr}_{\mathbb{R}}\circ\chi^{-1}$ satisfies $Z(H)=1$\;;\\
        \item given $H\in\mathcal{C}^{\infty}(M(\pi),\,\mathbb{R})$ such that $Z(H)=1$, then there exists an unique Hamiltonian section $h$ such that $H=\text{pr}_{\mathbb{R}}\circ\chi^{-1}$.
    \end{enumerate}
\end{prop}

\begin{proof}
    \begin{enumerate}[label=(\roman*)]
        \item Let $\eta\in M(\pi)$, so 
        \begin{eqnarray*}
            Z_{\eta}\big(\text{pr}_{\mathbb{R}}\circ\chi^{-1}\big) & = & \frac{d}{dt}\Big|_{t=0}\Big(\text{pr}_{\mathbb{R}}\circ\chi^{-1}\Big)\Big(\eta+t\big(\pi^*{vol}^{\Sigma}\big)_{\kappa\circ\rho(\eta)}\Big)\\
            & = & \frac{d}{dt}\Big|_{t=0} \text{pr}_{\mathbb{R}}\Bigg(\rho(\eta),\,\frac{\eta-h\circ\rho(\eta)+t\big(\pi^*{vol}^{\Sigma}\big)_{\kappa\circ\rho(\eta)}}{(\pi^*{vol}^{\Sigma})_{\kappa\circ\rho(\eta)}}\Bigg)\\
            & = & \frac{d}{dt}\Big|_{t=0}\Bigg(\frac{\eta-h\circ\rho(\eta)}{(\pi^*{vol}^{\Sigma})_{\kappa\circ\rho(\eta)}}+t\Bigg)=1.
        \end{eqnarray*}
        \item Let $H\in \mathcal{C}^{\infty}(M(\pi),\,\mathbb{R})$ such that for all $\eta\in M(\pi)$, $Z_{\eta}(H)=1$. Thus for $\eta\in M(\pi)$, we have $dH_{\eta}(Z)=Z_{\eta}(H)=1$, so the subspace $H^{-1}(\{0\})=\{H=0\}$ is a submanifold of codimension $1$. We show that $\rho\big|_{\{H=0\}}$ is a bijection, starting with injectivity.\\
        Let $\eta,\,\lambda\in \{H=0\}$ such that $\rho\big|_{\{H=0\}}(\eta)=\rho\big|_{\{H=0\}}(\lambda)=z\in P(\pi)$. Then there exists $u\in \mathbb{R}$ such that $\eta=\lambda+u(\pi^*{vol}^{\Sigma})_{\kappa\circ\rho(\eta)}$. But, by Lemma\nolinebreak[4] \ref{Flot Z},
        \begin{eqnarray*}
            0 & = & H(\eta)\\
            & = & H\big(\lambda+u(\pi^*{vol}^{\Sigma})_{\kappa\circ\rho(\eta)}\big)\\
        & = & H(\lambda)+u.
        \end{eqnarray*} 
        But $\lambda\in\{H=0\}$, so $u=0$ and $\eta=\lambda$.\\
        Let us now show surjectivity.\\
        Let $z\in P(\pi)$, then there exist $\eta\in M(\pi)$ such that $\rho(\eta)=z$ and $c\in \mathbb{R}$ such that $H(\eta)=c$. Thus 
        \begin{eqnarray*}
            H\big(\eta-c(\pi^*{vol}^{\Sigma})_{\kappa\circ\rho(\eta)}\big)=H(\eta)-c=0.
        \end{eqnarray*}
        So, $\eta-c(\pi^*{vol}^{\Sigma})_{\kappa\circ\rho(\eta)}\in \{H=0\}$ and $\rho\big|_{\{H=0\}}\big(\eta-c(\pi^*{vol}^{\Sigma})_{\kappa\circ\rho(\eta)}\big)=\rho(\eta)=z$, thus $\rho\big|_{\{H=0\}}$ is surjective and therefore a bijection.\\
        Furthermore, we have for $\eta\in M(\pi)$, $\big(T_{\eta}\rho\big)\big(T_{\eta}M(\pi)\big)=T_{\rho(\eta)}P(\pi)$ and $T_{\eta}M(\pi)=T_{\eta}\{H=0\}\oplus\big<Z_{\eta}\big>_{\mathbb{R}}$. Since $\big(T_{\eta}\rho\big)(Z_{\eta})=0$, we conclude that $\rho\big|_{\{H=0\}}$ is a diffeomorphism.\\
        We define $h=\Big(\rho\big|_{\{H=0\}}\Big)^{-1}$ and show its unicity. \raggedright If $\Tilde{h}$  is a Hamiltonian section that defines the function $\Tilde{H}$ and $\Tilde{H}=H$, one has $\text{im}(\Tilde{h})=\{\Tilde{H}=0\}=\{H=0\}=\text{im}(h)$ and thus the sections $\Tilde{h}$ and $h$ coincide.
    \end{enumerate}
\end{proof}

\section{Hamiltonian formulations of classical field theory}\label{Section : Hamiltonian formulations of Classical Field Theory}

In this section we first introduce the canonical forms on multicotangent bundles and certain notions central to the different formulations of \textit{Hamilton's equation} for classical fields, i.e. sections $\psi$ resp. $\widetilde{\psi}$ of the bundles $M(\pi)\rightarrow \Sigma$ resp. $P(\pi)\rightarrow \Sigma$.\\
Then we show in a completely self-contained way the equivalence of five different approaches to define a solution of a classical field theory problem\:: satisfaction of the equation $\widetilde{\psi}^*(\iota_X\omega_h)=0$ for all vector fields $X$ in $P(\pi)$ resp. of the (local) Hamilton-Volterra equations in standard coordinates resp. $\psi$ being a vortex $n$-plane resp. $\widetilde{\psi}$ satisfying Hamilton's principle of least action (\textit{extremal section}) resp. finally, $\psi$ satisfying the dynamical Hamilton-de Donder-Weyl equation (given a volume form on the source $\Sigma$ of $\psi$).\\

\noindent We conclude by giving a meaning to the dynamical Hamilton-de Donder-Weyl equation for a map from a $k$-dimensional source manifold $\Sigma$ to a $n$-plectic manifold $M$, given a differential form of degree $n$-$k$ on $M$ and a co-volume on $\Sigma$, greatly generalizing the case $M=M(\pi)$, $k=n$ considered hitherto in this article (i.e. the case of classical field theories with field values in multicotangent bundles).\\
Let us note that part of these equivalences are folklore in the field but with proofs scattered through the literature and sometimes lacking a completely rigorous formulation. (Important contributions to these equivalences are notably made in \cite{Carinena&Crampin&Ibort}, \cite{Echeverria&Leon&MunozLecanda&RomanRoy}, \cite{Goldschmidt&Sternberg}, \cite{Helein} and subsumed in \cite{Ryvkin&Wurzbacher} without complete proofs.)
\\\\
In this section, we are in the situation of Diagram \ref{diagram} with a standard coordinate system and fixed Hamiltonian section $h$.

\begin{defn}
    \raggedright On $M(\pi)$ we call \textit{Liouville form} the tautological form $\Theta$, i.e. for all $\eta\in M(\pi)$ over $y\in E$, 
    \begin{eqnarray*}
        \Theta_{\eta}(u_1,\,\ldots,\,u_n):= \eta_y\big((\kappa\circ \rho)_{*_y}u_1,\,\ldots,\,(\kappa\circ \rho)_{*_y}u_n\big).
    \end{eqnarray*}
    In coordinates we have $\Theta=pd^nx+\sum\limits_{a=1}^N\sum\limits_{\mu=1}^np^{\mu}_adq^a\wedge d^{n-1}\widehat{x^{\mu}}$.
\end{defn}

\begin{rem}
    \begin{enumerate}[label=(\roman*)]
        \item The Liouville form satisfies, for all $n$-forms $\alpha$ on $E$ \big(having values in $M(\pi)$\big), $\alpha^*\Theta=\alpha$.
        \item The negative of the differential of the Liouville form will be noted $\omega=\nolinebreak[4]-d\Theta$.
        \item If $n=1$, we are in the case of time-dependent classical mechanics with $\Theta$ the classical Liouville form and $\omega$ the symplectic form on $\Lambda^1_2T^*E=\nolinebreak[4]T^*E$. \Big(Note that locally $\omega=dx\wedge dp+\sum\limits_{a=1}^Ndq^a\wedge dp_a$\Big).
        \item We will call $\Theta_h=h^*\Theta$ and $\omega_h=h^*\omega$ the ensuing pull-backs. In a standard coordinate chart, we have
        \begin{eqnarray*}
            \omega_h=h^*\omega=d\mathcal{H}\wedge d^nx-\sum_{\eta=1}^n\sum_{a=1}^Ndp^{\eta}_a\wedge dq^a\wedge d^{n-1}\widehat{x^{\eta}}.
        \end{eqnarray*}
    \end{enumerate}
\end{rem}

\begin{defn}
    Given a coordinate patch $U\subset \Sigma$ and setting ${vol}^{U}:=d^nx$, we define $\mathcal{H}=\mathcal{H}^U$ on $\tau^{-1}(U)$ as in Section \ref{Section : Hamiltonian sections versus Hamiltonian functions}. We named this a local Hamilton-Volterra function. Furthermore, we call the following equations for a section $\widetilde{\psi}$ of $\tau$ over $U$ the \textit{(local) Hamilton-Volterra equations} :
    \begin{eqnarray*}
        -\frac{\partial \mathcal{H}}{\partial q^a}\big(\widetilde{\psi}(x)\big)  & = &  \sum_{\mu=1}^n\frac{\partial (p^\mu_a\circ \widetilde{\psi})}{\partial x^\mu}(x)\;;\\
        \frac{\partial \mathcal{H}}{\partial p^\mu_a}\big(\widetilde{\psi}(x)\big) & = & \frac{\partial (q^a\circ \widetilde{\psi})}{\partial x^\mu}(x).
    \end{eqnarray*}
\end{defn}

\begin{rem}
    The above equations appear apparently for the first time in the paper \cite{Volterra1890EQ} from Volterra.
\end{rem}

\begin{rem}
    We will in the sequel always denote a section of $\tau$ over $U$ by $\widetilde{\psi}$ and by $\psi=h\circ \widetilde{\psi}$ the induced section of $(\tau\circ\rho)$ over $U$.
\end{rem}

\begin{defn}
    Let $U\subset \Sigma$ be open and $\Psi$ a section of $(\tau\circ\rho)$. We call $\Psi$ a \textit{(local) vortex $n$-plane} if 
    \begin{enumerate}[label=(\roman*)]
        \item $\Psi(U)\subset \text{im}(h)$ and
        \item $\forall x\in U,\,\forall \gamma_x\in \Lambda^nT_x\Sigma$ one has $\iota_{(\Psi_*)_x(\gamma_x)}\omega_{\Psi(x)}=0$ as a functional on $T_{\Psi(x)}\text{im}(h)$.
    \end{enumerate}
\end{defn}

\begin{rem}
    In \cite{Kijowski}, Kijowski calls a vortex $n$-plane a \textit{state}.
\end{rem}

\begin{thm}\label{Equivalence Theorem}
\raggedright Let $\widetilde{\psi}$ be a section of $\tau$ defined on an open subset $U\subset \Sigma$ and $\psi=h\circ\widetilde{\psi}$, then the following are equivalent\::
    \begin{enumerate}[label=(\roman*)]
        \item $\widetilde{\psi}^*(\iota_X\omega_h)=0$ for all $X$ $\tau$-vertical in $\mathfrak{X}\big(P(\pi)\big)$ (i.e. $X$ is in the kernel of $\tau_*$)\;;
        \item $\widetilde{\psi}^*(\iota_X\omega_h)=0$ for all $X\in \mathfrak{X}\big(P(\pi)\big)$\;;
        \item in standard coordinates $(x^{\mu},\,q^a,\,p^{\mu}_a)$ on $P(\pi)$, $\widetilde{\psi}$ satisfies the Hamilton-Volterra equations\;;
        \item the section $\psi$ is a vortex $n$-plane.
    \end{enumerate}
\end{thm}

\begin{rem}
    The conditions $(i)$ and $(ii)$ can equivalently be formulated using only vector fields with compact support since the conditions are local.
\end{rem}

\begin{proof}
    \raggedright $(i)\implies (ii)$ If $p\in \text{im}(\widetilde{\psi})$, so $T_pP(\pi)=V_p(\tau)\oplus T_p\big(\text{im}(\widetilde{\psi})\big)$. Then, for $X\in \mathfrak{X}\big(P(\pi)\big)$, we have the following decomposition 
    \begin{eqnarray*}
        X_p=\big(X_p-T_p(\widetilde{\psi}\circ\tau)(X_p)\big)+T_p(\widetilde{\psi}\circ\tau)(X_p)=X_p^V+X_p^{\widetilde{\psi}}
    \end{eqnarray*}
    where $X_p^V\in V_p(\tau)$ is the vertical part of the vector field and $X_p^{\widetilde{\psi}}\in T_p\big(\text{im}(\widetilde{\psi})\big)$.\\
    Moreover 
    \begin{eqnarray*}
        \widetilde{\psi}^*(\iota_X\omega_h)=\widetilde{\psi}^*(\iota_{X^V}\omega_h)+\widetilde{\psi}^*(\iota_{X^{\widetilde{\psi}}}\omega_h)=\widetilde{\psi}^*(\iota_{X^{\widetilde{\psi}}}\omega_h)=0
    \end{eqnarray*}
    because $X_p^{\widetilde{\psi}}\in T_p\big(\text{im}(\widetilde{\psi})\big)$ and $\dim\big(\text{im}(\widetilde{\psi})\big)=n$. The converse is obvious.\\

    $(ii) \iff (iii)$ Write, in local coordinates 
    \begin{eqnarray*}
        X=\sum_{\mu=1}^n \alpha^{\mu}\frac{\partial}{\partial x^{\mu}} +\sum_{a=1}^N\beta^a\frac{\partial}{\partial q^a}+\sum_{\mu=1}^n\sum_{a=1}^N\gamma^{\mu}_a\frac{\partial}{\partial p^{\mu}_a}\in\mathfrak{X}\big(P(\pi)\big).
    \end{eqnarray*}
    Then
    \begin{eqnarray*}
        \iota_X\omega_h = & &  -\sum_{a=1}^N\sum_{\eta,\nu=1}^n\alpha^{\eta}\bigg(dp_a^{\nu}\wedge dq^a\wedge d^{n-2}\widehat{x^{\eta,\nu}}\bigg)\\
        & - & \sum_{a=1}^N\sum_{\eta=1}^n\alpha^{\eta}\bigg(\sum_{\nu=1}^n\frac{\partial \mathcal{H}}{\partial p^{\nu}_a} dp^{\nu}_a\wedge d^{n-1}\widehat{x^{\eta}}+\frac{\partial \mathcal{H}}{\partial q^a}dq^a\wedge d^{n-1}\widehat{x^{\eta}}\bigg)\\
        & + & \sum_{\nu=1}^n\sum_{a=1}^N\beta^adp^{\nu}_a\wedge d^{n-1}\widehat{x^{\nu}}+\sum_{a=1}^N\beta^a\frac{\partial \mathcal{H}}{\partial q^a}d^nx\\
        & - & \sum_{\nu=1}^n\sum_{a=1}^N\gamma^{\nu}_a\bigg(dq^a\wedge d^{n-1}\widehat{x^{\nu}}-\frac{\partial \mathcal{H}}{\partial p^{\nu}_a}d^nx\bigg)
    \end{eqnarray*}
    where $d^{n-2}\widehat{x^{\eta,\nu}}=\iota_{\frac{\partial}{\partial x^{\eta}}}\iota_{\frac{\partial}{\partial x^{\nu}}}d^nx$.
    If $\widetilde{\psi}=\big(x^{\nu},\,q^a(x^{\eta}),\,p_a^{\nu}(x^{\eta})\big)$, the first term of $\widetilde{\psi}^*(\iota_X\omega_h)$ is\::
    \begin{eqnarray*}
        & & -\widetilde{\psi}^*\Bigg(\sum_{a=1}^N\sum_{\eta,\nu=1}^n  \alpha^{\eta}\bigg(dp_a^{\nu}\wedge dq^a\wedge d^{n-2}\widehat{x^{\eta,\nu}}\bigg)\Bigg)\\
        & = & \widetilde{\psi}^*\Bigg(\sum_{a=1}^N\sum_{\eta,\nu =1,\eta\neq \nu}^n\alpha^{\eta}\bigg(dp_a^{\nu}\wedge dq^a\wedge d^{n-2}\widehat{x^{\eta,\nu}}\bigg)\Bigg)\\
        & = & \sum_{a=1}^N\sum_{\eta,\nu=1,\,\eta\neq \nu}^n\alpha^{\eta}\bigg(\frac{\partial (q^a\circ\widetilde{\psi})}{\partial x^{\nu}} \frac{\partial (p^{\nu}_a\circ\widetilde{\psi})}{\partial x^{\eta}}-\frac{\partial (q^a\circ\widetilde{\psi})}{\partial x^{\eta}} \frac{\partial (p^{\nu}_a\circ\widetilde{\psi})}{\partial x^{\nu}}\bigg) d^nx\\
        & = & \sum_{a=1}^N\sum_{\eta,\nu=1}^n\alpha^{\eta}\bigg(\frac{\partial (q^a\circ\widetilde{\psi})}{\partial x^{\nu}} \frac{\partial (p^{\nu}_a\circ\widetilde{\psi})}{\partial x^{\eta}}-\frac{\partial (q^a\circ\widetilde{\psi})}{\partial x^{\eta}} \frac{\partial (p^{\nu}_a\circ\widetilde{\psi})}{\partial x^{\nu}}\bigg) d^nx.
    \end{eqnarray*}
    We calculte similarly the other terms. Finally, we obtain $\widetilde{\psi}^*(\iota_X\omega_h)=$
    \begin{eqnarray*}
        & & \sum_{a=1}^N\sum_{\eta,\nu=1}^n\alpha^{\eta}\bigg(\frac{\partial (q^a\circ\widetilde{\psi})}{\partial x^{\nu}} \frac{\partial (p^{\nu}_a\circ\widetilde{\psi})}{\partial x^{\eta}}-\frac{\partial (q^a\circ\widetilde{\psi})}{\partial x^{\eta}} \frac{\partial (p^{\nu}_a\circ\widetilde{\psi})}{\partial x^{\nu}}\bigg) d^nx\\
        & & - \sum_{a=1}^N\sum_{\eta=1}^n\alpha^{\eta}\bigg(\sum_{\nu=1}^n\frac{\partial \mathcal{H}}{\partial p^{\nu}_a}\Big|_{\widetilde{\psi}} \frac{\partial (p^{\nu}_a\circ\widetilde{\psi})}{\partial x^{\eta}}d^nx+\frac{\partial \mathcal{H}}{\partial q^a}\Big|_{\widetilde{\psi}}\frac{\partial (q^a\circ \widetilde{\psi})}{\partial x^{\eta}} d^nx\bigg)\\
        & & +  \sum_{\nu=1}^n\sum_{a=1}^N\beta^a\frac{\partial (p^{\nu}_a\circ\widetilde{\psi})}{\partial x^{\nu}} d^nx+\sum_{a=1}^N\beta^a\frac{\partial \mathcal{H}}{\partial q^a}\Big|_{\widetilde{\psi}}d^nx\\
        & & - \sum_{\nu=1}^n\sum_{a=1}^N\gamma^{\nu}_a\bigg(\frac{\partial (q^a\circ\widetilde{\psi})}{\partial x^{\nu}} d^nx-\frac{\partial \mathcal{H}}{\partial p^{\nu}_a}\Big|_{\widetilde{\psi}}d^nx\bigg),
    \end{eqnarray*}
    thus
    \begin{eqnarray*}
        \widetilde{\psi}^*(\iota_X\omega_h) & = & \sum_{a=1}^N\sum_{\eta,\nu=1}^n\alpha^{\eta}\bigg(\frac{\partial (q^a\circ \widetilde{\psi})}{\partial x^{\nu}}-\frac{\partial \mathcal{H}}{\partial p_a^\nu}\Big|_{\widetilde{\psi}}\bigg)\frac{\partial (p^{\nu}_a\circ\widetilde{\psi})}{\partial x^{\eta}}d^nx\\
        & & - \sum_{a=1}^N\sum_{\eta=1}^n\alpha^{\eta}\bigg(\sum_{\nu=1}^n\frac{\partial (p_a^{\nu}\circ \widetilde{\psi})}{\partial x^{\nu}}+\frac{\partial \mathcal{H}}{\partial q^a}\Big|_{\widetilde{\psi}}\bigg)\frac{\partial (q^a\circ \widetilde{\psi})}{\partial x^{\eta}}d^nx\\
        & & + \sum_{a=1}^N\Bigg(\beta^a\bigg(\sum_{\nu=1}^n\frac{\partial (p_a^{\nu}\circ \widetilde{\psi})}{\partial x^{\nu}}+\frac{\partial \mathcal{H}}{\partial q^a}\Big|_{\widetilde{\psi}}\bigg)d^nx\\
        & & - \sum_{\nu=1}^n\gamma_a^{\nu}\bigg(\frac{\partial (q^a\circ \widetilde{\psi})}{\partial x^{\nu}}-\frac{\partial \mathcal{H}}{\partial p_a^{\nu}}\Big|_{\widetilde{\psi}}\bigg)d^nx\Bigg).
    \end{eqnarray*}
    We conlude that the condition $\widetilde{\psi}^*(\iota_X\omega_h)=0$ for all $X\in \mathfrak{X}\big(P(\pi)\big)$ is equivalent to the Hamilton-Volterra equations.\\

    $(ii)\iff (iv)$ For $X\in\mathfrak{X}\big(P(\pi)\big),\, x\in U$ and $\gamma_x\in \Lambda^nT_x\Sigma$ we have
    \begin{eqnarray*}
        \big(\widetilde{\psi}^*(\iota_X\omega_h)\big)_x(\gamma_x) & = & (\omega_h)_{\widetilde{\psi}(x)}\big(X_{\widetilde{\psi}(x)},(\widetilde{\psi}_*)_x(\gamma_x)\big) \\
        & = & (-1)^n\big(\iota_{(\widetilde{\psi}_*)_x(\gamma_x)}(\omega_h)_{\widetilde{\psi}(x)}\big)(X_{\widetilde{\psi}(x)}).
    \end{eqnarray*}
    Thus $\widetilde{\psi}^*(\iota_X\omega_h)=0$ for all $X\in \mathfrak{X}\big(P(\pi)\big)$ is equivalent to $\big(\iota_{(\widetilde{\psi}_*)_x(\gamma_x)}(\omega_h)_{\widetilde{\psi}(x)}\big)(v_{\widetilde{\psi}(x)})=0$, for all $x\in U\subset \Sigma$, for all $\gamma_x\in \Lambda^nT_x\Sigma$ and for all $v_{\widetilde{\psi}(x)}\in T_{\widetilde{\psi}(x)}P(\pi)$.
\end{proof}

\noindent In the sequel, the notation $V\subset \subset U$ will be used if $V$ is a relatively compact and open subset of $U$.

\begin{defn}
    Let $U\subset \Sigma$ be an orientable open subset. For a section $\widetilde{\psi}$ of $\tau\:: P(\pi)\big|_U\rightarrow U$, we define, upon fixing an orientation on $U$ for all $V\subset \subset U$ open with smooth boundary,
    \begin{eqnarray*}
         \mathbb{H}_V[\widetilde{\psi}]=\int_V\widetilde{\psi}^*\Theta_h=\int_V\psi^*\Theta.
    \end{eqnarray*}
    We call $\widetilde{\psi}$ \textit{extremal} if $\widetilde{\psi}$ is critical for $\mathbb{H}_V$ for all such $V$ in the following sense\::
    for all smooth families of sections $(\widetilde{\psi}_t)_{|t|< \varepsilon}$ with $\varepsilon>0$ and $\widetilde{\psi}_0=\widetilde{\psi}$ such that $\forall t$, $\widetilde{\psi}_t=\widetilde{\psi}$ outside a fixed compact subset of $V$, 
    \begin{eqnarray*}
        \frac{d}{dt}\Big|_{t=0}\mathbb{H}_V[\widetilde{\psi}_t]=0.
    \end{eqnarray*}
\end{defn}

\begin{rem}
    We can drop the assumption that the boundary of $V$ is smooth since there always exists a $W\subset \subset U$ with smooth boundary containing a given relatively compact $V\subset \subset U$.
\end{rem}

\begin{rem}
    Smooth families as in the preceeding definition are called \textit{admissible families (relative to $V\subset \subset U)$} in the sequel. Note that even if $(\widetilde{\psi}_t)_{|t|< \varepsilon}$ are initially only defined over $V$ they all can be extended to $U$ upon setting $\widetilde{\psi}_t(x)=\widetilde{\psi}(x)$ for all $x\in U\backslash V$.
\end{rem}

\begin{ex}
    \begin{enumerate}[label=(\roman*)]
        \item Let us check that extremality already implies  the Hamilton equations in classical mechanics. Let $Q$ be open in $\mathbb{R}^N$, $I$ open in $\mathbb{R}$ and $\mathcal{H}\::P(\pi)=I\times T^*Q\rightarrow \mathbb{R}$ a time-dependant Hamilton function. Then $(x,\,q^a,\,p_a)$ are global coordinates on $P(\pi)$ and $\widetilde{\psi}(x)=\big(x,\,q^a(x),\,p_a(x)\big)$ and $h(x,\,q^a,\,p_a)=\big(x,\,q^a,\,p_a,\,-\mathcal{H}(x,\,q^a,\,p_a)\big)$, \raggedright implying that $\psi^*\Theta=\psi^*(pdx+p_adq^a)=-\mathcal{H}dx+p_a\frac{\partial q^a}{\partial x}dx$. Thus, we obtain
        \begin{eqnarray*}
            \mathbb{H}_V[\widetilde{\psi}]=\int_V\Big(p_a\frac{\partial q^a}{\partial x}-\mathcal{H}\Big)dx
        \end{eqnarray*} for an admissible smooth family of sections $(\widetilde{\psi}_t)_{|t|< \varepsilon}$ such that $\widetilde{\psi}_t=\widetilde{\psi}+t\widetilde{\varphi}$ on $V\subset Q$ where $\widetilde{\varphi}(x)=\big(x,\,u^a(x),\,v_a(x)\big)$,
        \begin{eqnarray*}
            \raggedright 0=\frac{d}{dt}\Big|_{t=0}\mathbb{H}_V[\widetilde{\psi}_t]= \int_Vu^a\Big(-\frac{\partial p_a}{\partial x}-\frac{\partial \mathcal{H}}{\partial q^a}\Big)dx +\int_Vv_a\Big(\frac{\partial q^a}{\partial x}-\frac{\partial \mathcal{H}}{\partial p_a}\Big)dx,
        \end{eqnarray*}
        i.e. Hamilton's equations for $\big(q^a(x),\,p_a(x)\big)$.
        \item Let us explain, by an example, why criticality has to be replaced by extremality. Let $\Sigma=\mathbb{R}^2$ and $E=\mathbb{R}^2\times \mathbb{R}\xrightarrow{\pi} \mathbb{R}^2$ with coordinates $(x^1,\,x^2,\,q,\,p^1,\,p^2)$ on $P(\pi)\cong \mathbb{R}^2\times \mathbb{R}\times \mathbb{R}^2$ and $\mathcal{H}=\frac{1}{2}\big((p^1)^2+(p^2)^2\big)$. Then the Hamilton-Volterra equations 
        \begin{eqnarray*}
            p^{\mu}=\frac{\partial q}{\partial x^{\mu}} & \text{and} & \frac{\partial^2q}{\partial x^1}+\frac{\partial^2q}{\partial x^2}=0
        \end{eqnarray*}
        \raggedright are easily obtained via the localized functionals $\mathbb{H}_V$, but there is no non-constant harmonic function $q$ on $\Sigma=\mathbb{R}^2$ such that for $\widetilde{\psi}(x)=\Big(x,\,q(x),\,\frac{\partial q}{\partial x^{\mu}}\Big)$,
        \begin{eqnarray*}
            \mathbb{H}_{\Sigma}[\widetilde{\psi}] & = & \int_{\mathbb{R}^2}\Big(-\mathcal{H}\big(\widetilde{\psi}(x)\big)+p^1\frac{\partial q}{\partial x^1}+p^2\frac{\partial q}{\partial x^2}\Big) dx^1\wedge dx^2\\
            & = & \frac{1}{2}\int_{\mathbb{R}^2}||\overrightarrow{\nabla}q||^2dx^1\wedge dx^2
        \end{eqnarray*}
        converges. Thus passing to extremality is necessary.
    \end{enumerate}
\end{ex}

\begin{rem}
    \begin{enumerate}[label=(\roman*)]
        \item Extremality is independant of the choice of an orientation on $U$.
        \item Being extremal in the above sense is a multidimensional analogue of satisfying Hamilton's principle of least action in classical mechanics, compare \cite{Arnold}, p.243 (see also the corresponding footnote on page 246).
    \end{enumerate}
\end{rem}

\begin{thm}
    Let $U\subset \Sigma$ be open and orientable and $\widetilde{\psi}$ a section of $\tau$ over $U$. Then the conditions of Theorem \ref{Equivalence Theorem} are also equivalent to $\widetilde{\psi}$ being an extremal section.
\end{thm}

\begin{proof}
    \raggedright We will show that $\widetilde{\psi}$ being an extremal section is equivalent to the first condition of Theorem \ref{Equivalence Theorem}\:: $\widetilde{\psi}^*(\iota_X\omega_h)=0$ for all $X$ $\tau$-vertical in $\mathfrak{X}\big(P(\pi)\big)$.\\
    Let first $X$ be a vertical vector field with compact support on $P(\pi)$. Take $V\subset \subset U$ open with smooth boundary such that $\text{supp}(X)\subset \tau^{-1}(V)$. We call $\sigma_t$ the flow of this vector field and set $\widetilde{\psi}_t=\sigma_t\circ\widetilde{\psi}$. Remark that $(\widetilde{\psi}_t)_{|t|<\varepsilon}$ is an admissible family of sections of $\tau$. Then
    \begin{eqnarray*}
        \raggedright \frac{d}{dt}\Big|_{t=0}\int_V\widetilde{\psi}_t^*\Theta_h & = & \frac{d}{dt}\Big|_{t=0}\int_V\widetilde{\psi}^*(\sigma_t^*\Theta_h)=\int_V\widetilde{\psi}^*\Big(\lim_{t\rightarrow 0}\frac{\sigma_t^*\Theta_h-\Theta_h}{t}\Big) \\
        & = & \int_V\widetilde{\psi}^*(\mathcal{L}_X\Theta_h)= \int_V\widetilde{\psi}^*(\iota_Xd\Theta_h+d\iota_X\Theta_h)\\
        & = & -\int_V\widetilde{\psi}^*(\iota_X\omega_h)+\int_Vd[\widetilde{\psi}^*(\iota_X\Theta_h)]\\
        & = & -\int_V\widetilde{\psi}^*(\iota_X\omega_h)+\int_{\partial V}\widetilde{\psi}^*(\iota_X\Theta_h)=-\int_V\widetilde{\psi}^*(\iota_X\omega_h).
    \end{eqnarray*}
    Note that the integral over $\partial V$ is zero since for $x\in\partial V$, $\widetilde{\psi}(x)$ cannot be in supp($X$) and thus the form $\widetilde{\psi}^*(\iota_X\Theta_h)$ vanishes everywhere on $\partial V$.\\
    Thus extremality of $\widetilde{\psi}$, i.e. the condition
    \begin{eqnarray*}
        \frac{d}{dt}\Big|_{t=0}\int_V\widetilde{\psi}_t^*\Theta_h=0, \text{ implies } \widetilde{\psi}^*(\iota_X\omega_h)=0.
    \end{eqnarray*}
    Let now condition $(i)$ of Theorem \ref{Equivalence Theorem} be satisfied, $V\subset\subset U$ and $(\widetilde{\psi}_t)_{|t|<\varepsilon}$ an admissible family of sections. We define the vertical vector field $\xi$ along $\widetilde{\psi}_0=\widetilde{\psi}$ by $\xi(x)=\frac{d}{dt}\big|_{t=0}\widetilde{\psi}_t(x)\in T_{\widetilde{\psi}(x)}P(\pi)$. We extend $\xi$ from $\widetilde{\psi}(U)$ to a vertical field $X$ on $P(\pi)$ with compact support supp($X$)$\subset \tau^{-1}(V)$.\\
    Let us note the following identity that is directly implied by formula $(22.1)$ in \cite{Guillemin}\::
    \begin{eqnarray*}
        \frac{d\widetilde{\psi}_t^*\Theta_h}{dt}\Big|_{t=0}=d\big(\widetilde{\psi}_0^*(\iota_{\xi}\Theta_h)\big)+\widetilde{\psi}_0^*(\iota_{\xi}d\Theta_h),
    \end{eqnarray*}
    where for $F\::M\rightarrow N$ smooth, $\alpha\in \Omega^{p+1}(N)$ and $\xi$ a vector field along $F$, setting
    \begin{eqnarray*}
        \big(F^*(\iota_{\xi}\alpha)\big)_x(u_1,\,\ldots,\,u_n):=\alpha_{F(x)}\big(\xi_x,\,(F_*)_x(u_1),\,\ldots,\,(F_*)_x(u_n)\big)
    \end{eqnarray*}
    for all $u_1,\,\ldots,\,u_n\in T_xM$ defines a $p$-form on $M$.\\
    Now we compute\::
    \begin{eqnarray*}
        \frac{d}{dt}\Big|_{t=0}\mathbb{H}_V[\widetilde{\psi}_t] & = & \int_V\frac{d}{dt}\Big|_{t=0}(\widetilde{\psi}_t^*\Theta_h)\\
        & = & \int_V\Big(d\big(\widetilde{\psi}^*_0(\iota_{\xi}\Theta_h)\big)+\widetilde{\psi}^*_0(\iota_{\xi}d\Theta_h)\Big)\\
        & = & \int_{\partial V}\widetilde{\psi}^*(\iota_{\xi}\Theta_h)-\int_V\widetilde{\psi}^*(\iota_{\xi}\omega_h)\\
        & = & -\int_V\widetilde{\psi}^*(\iota_X\omega_h).
    \end{eqnarray*}
    Here the boundary integral vanishes since for $x\in\partial V,\,\xi(x)=0$ and in the last equation we can replace $\xi$ by $X$ since for $x\in V$ one has $X\big(\widetilde{\psi}(x)\big)\!=\!\xi(x)$.\\
    The assumption $\widetilde{\psi}^*(\iota_X\omega_h)=0$ now implies that for every admissible family $(\widetilde{\psi}_t)_{|t|<\varepsilon}$, one has $\frac{d}{dt}\big|_{t=0}\mathbb{H}_V[\widetilde{\psi}_t]=0$.
\end{proof}

\begin{thm}\label{section dynamical HDW}
    Let ${vol}^{\Sigma}$ be a volume form on $\Sigma$ and $H$ the Hamiltonian function associated to a Hamiltonian section $h$ and ${vol}^{\Sigma}$. Denote by $\gamma\in\nolinebreak[4]\mathfrak{X}^n(\Sigma)$ the $n$-vector field dual to ${vol}^{\Sigma}$, i.e. $({vol}^{\Sigma})_x(\gamma_x)=1,\,\forall x\in \Sigma$. Let, furthermore, $U$ be open in $\Sigma$ and $\widetilde{\psi}$ a section of $\tau$ over $U$. Then the conditions of Theorem \ref{Equivalence Theorem} are equivalent to 
    \begin{eqnarray*}
        \forall x\in U,\, \iota_{(\psi_*)_x(\gamma_x)}\omega_{\psi(x)}=(-1)^{n+1}(dH)_{\psi(x)}.
    \end{eqnarray*}
    (Again, $\psi=h\circ \widetilde{\psi}$ and $\widetilde{\psi}=\rho\circ\psi$ here.)
\end{thm}

\begin{proof}
    We will show that the section $\psi$ being a vortex $n$-plane is equivalent to satisfying the equation 
    \begin{eqnarray*}
        \forall x\in U,\, \iota_{(\psi_*)_x(\gamma_x)}\omega_{\psi(x)}=(-1)^{n+1}(dH)_{\psi(x)}.
    \end{eqnarray*}
    Let $W=\{H=0\}\subset U$. On $T_{\psi(x)}W$, we have $\iota_{(\psi_*)_x(\gamma_x)}\omega_{\psi(x)}=0$. As a functional 
    on $T_{\psi(x)}M(\pi),\,\iota_{(\psi_*)_x(\gamma_x)}\omega_{\psi(x)}=g(x)(dH)_{\psi(x)},\,\forall x\in U,$ where $g:U\rightarrow \mathbb{R}$ is a smooth function. We know by Proposition \ref{Proposition H=h+p} that $H=\mathcal{H}+p$ for a certain function $\mathcal{H}$, so we have $\iota_{\big(\frac{\partial}{\partial p}\big)_{\psi(x)}}g(x)(dH)_{\psi(x)}=g(x)$ and thus
    \begin{eqnarray*}
        1 & = & \iota_{(\psi_*)_x(\gamma_x)}d^nx_{\psi(x)} = -\iota_{(\psi_*)_x(\gamma_x)}\iota_{\frac{\partial}{\partial p}}\omega_{\psi(x)}\\
        & = & (-1)^{n+1}\iota_{\frac{\partial}{\partial p}}\iota_{(\psi_*)_x(\gamma_x)}\omega_{\psi(x)}\\
        & = & (-1)^{n+1} \iota_{\frac{\partial}{\partial p}}\big(g(x)(dH)_{\psi(x)}\big)\\
        & = & (-1)^{n+1}g(x).
    \end{eqnarray*}
    We conclude that $g(x)=(-1)^{n+1}$ and that the announced equation holds true.\\
    For the converse, remark that $\psi(U)\subset \text{im}(h) = \{H=0\}$ and we find the equation for a vortex $n$-plane.
\end{proof}

\begin{rem}
In standard coordinates the preceding theorem gives back, of course,  the Hamilton-Volterra  equations plus the following \textit{energy equations} for $\mu = 1,\ldots,n$ 

\begin{eqnarray*}
    \frac{\partial p}{\partial x^{\mu}}=-\frac{\partial H}{\partial x^{\mu}}-\sum_{a,\mu\neq \nu} \Big(\frac{\partial q^a}{\partial x^{\nu}}\frac{\partial p^{\nu}_a}{\partial x^{\mu}}-
    \frac{\partial q^a}{\partial x^{\mu}}\frac{\partial p^{\nu}_a}{\partial x^{\nu}}\Big).
\end{eqnarray*}
These equations are already implied by the Hamilton-Volterra equations and thus redundant. In a closely related setting, this was already observed in \cite{Kijowski} on p. 108.
\end{rem}

\noindent We close this article by generalizing the condition formulated in the preceding theorem to a very large class of manifolds and maps.

\begin{defn}\label{Final Definition}
    \begin{enumerate}[label=(\roman*)]
        \raggedright\item A \textit{multisymplectic manifold} $(M,\,\omega)$ is a pair consisting of a manifold  $M$  and a closed differential form 
        $\omega\in \Omega_{cl}^{n+1}(M)$ with  $n\geq 1$ satisfying the following non-degeneracy condition\:: the map
        \begin{eqnarray*}
            \omega^{\flat}\:: TM & \rightarrow & \Lambda^nT^*M\\
            v & \mapsto & \omega^{\flat}(v)=\iota_v\omega
        \end{eqnarray*}
        is injective. For a fixed degree $n{+}1$ such a pair is also called an \textit{$n$-plectic manifold}. The form $\omega$ is sometimes called a \textit{multisymplectic form} or a \textit{multisymplectic structure} on $M$.
        \item Let $(M,\,\omega)$ be a $n$-plectic manifold. A couple $(H,\,X)\in \Omega^{n-k}(M)\times \mathfrak{X}^k(M)$ with 
        $1\leq k\leq n$ is a solution of the \textit{Hamilton-de Donder-Weyl equation (or HDW equation)} if 
        \begin{eqnarray*}
            \iota_X\omega=(-1)^{n+1-k}dH.
        \end{eqnarray*}
        \item Let $(M,\,\omega)$ be a $n$-plectic manifold, $\Sigma$ a $k$-dimensional manifold with $1\leq k\leq n$ and $H\in \Omega^{n-k}(M)$. A couple $(\gamma,\,\psi)$ whith $\gamma\in \mathfrak{X}^k(\Sigma)$ a co-volume and $\psi\::\Sigma\rightarrow M$ is called a solution of the \textit{dynamical Hamilton-de Donder-Weyl equation (or dynamical HDW equation)} if $\forall x\in \Sigma$
        \begin{eqnarray*}
            \iota_{(\psi_*)_x(\gamma_x)}\omega_{\psi(x)}=(-1)^{n+1-k}(dH)_{\psi(x)}.
        \end{eqnarray*}
    \end{enumerate}
\end{defn}

\begin{rem}
    \begin{enumerate}[label=(\roman*)]
        \item We learned the notion of HDW equation from \cite{Schreiber}, Definition 1.2.291.
        \item A solution of the dynamical HDW equation for $k=n$ on an $n$-plectic manifold is called a \textit{Hamiltonian $n$-curve} in \cite{Helein}.
        \item The manifold $M(\pi)$ is $n$-plectic and Theorem \ref{section dynamical HDW} can be stated by saying that a section $\psi$ satisfies the conditions of Theorem \ref{Equivalence Theorem} if and only if it satisfies the dynamical HDW equation for $H\in\nolinebreak[4]\Omega^0\big(M(\pi)\big)$.
        \item If $n\geq 2$, $\omega_h$ is $n$-plectic but note that this form mixes the geometry with a fixed dynamical system contrary to the symplectic case or the situation on $M(\pi)$.
        \item Some care is needed in aiming for a solution theory for dynamical HDW equations\:: if $H$ is closed, any map $\psi$ of rank constantly below $k$ solves the equation. The most natural non-degeneration condition seems to be to ask that $dH$ is (generically) non-vanishing. Note that there is another degeneration problem already present for the HDW equation\:: since, in general, there are non-vanishing $k$-vector fields $X$ with $\iota_X\omega=0$ for $k\geq2$ the ensuing couples $(0,\,X)$ solve the HDW equation. If $(\psi_*)_x(\gamma_x)=X\big(\psi(x)\big)$ for all $x\in\Sigma$ this yields a solution of the dynamical HDW equation, possibly with rank $(\psi_*)_x=k$ for all $x\in\Sigma$.
    \end{enumerate}
\end{rem}

\section*{Acknowledgment} Both authors profited from giving related talks in several seminars and conferences notably in Bia\l ystok, Lyon, Madrid, Metz and Paris. 

\bibliographystyle{spmpsci}
\bibliography{WW-References}

\end{document}